\documentclass[a4paper,12pt]{amsart}
\usepackage{amssymb,amsthm}
\usepackage[all]{xy}

\begin{document}

\newcommand{\opp}{\bowtie }
\newcommand{\pos}{\text {\rm pos}}
\newcommand{\supp}{\text {\rm supp}}
\newcommand{\End}{\text {\rm End}}
\newcommand{\diag}{\text {\rm diag}}
\newcommand{\Lie}{\text {\rm Lie}}
\newcommand{\Ad}{\text {\rm Ad}}
\newcommand{\car}{\mathcal R}
\newcommand{\Tr}{\rm Tr}
\newcommand{\Spec}{\text{\rm Spec}}

\def\ge{\geqslant}
\def\le{\leqslant}
\def\a{\alpha}
\def\b{\beta}
\def\c{\chi}
\def\g{\gamma}
\def \fg{{\mathfrak g}}
\def\G{\Gamma}
\def\d{\delta}
\def\D{\Delta}
\def\L{\Lambda}
\def\e{\epsilon}
\def\et{\eta}
\def\io{\iota}
\def\o{\omega}
\def\p{\pi}
\def\ph{\phi}
\def\ps{\psi}
\def\r{\rho}
\def\s{\sigma}
\def\t{\tau}
\def\th{\theta}
\def\k{\kappa}
\def\l{\lambda}
\def\z{\zeta}
\def\v{\vartheta}
\def\x{\xi}
\def\i{^{-1}}

\def\mapright#1{\smash{\mathop{\longrightarrow}\limits^{#1}}}
\def\mapleft#1{\smash{\mathop{\longleftarrow}\limits^{#1}}}
\def\mapdown#1{\Big\downarrow\rlap{$\vcenter{\hbox{$\scriptstyle#1$}}$}}

\def\ca{\mathcal A}
\def\cb{\mathcal B}
\def\cc{\mathcal C}
\def\cd{\mathcal D}
\def\ce{\mathcal E}
\def\cf{\mathcal F}
\def\cg{\mathcal G}
\def\ch{\mathcal H}
\def\ci{\mathcal I}
\def\cj{\mathcal J}
\def\ck{\mathcal K}
\def\cl{\mathcal L}
\def\cm{\mathcal M}
\def\cn{\mathcal N}
\def\co{\mathcal O}
\def\cp{\mathcal P}
\def\cq{\mathcal Q}
\def\car{\mathcal R}
\def\cs{\mathcal S}
\def\ct{\mathcal T}
\def\cu{\mathcal U}
\def\cv{\mathcal V}
\def\cw{\mathcal W}
\def\cz{\mathcal Z}
\def\cx{\mathcal X}
\def\cy{\mathcal Y}

\def\tz{\tilde Z}
\def\ty{\tilde Y}
\def\tl{\tilde L}
\def\tc{\tilde C}
\def\ta{\tilde A}
\def\tx{\tilde X}

\def \us {\underline{*}}
\def \os {\overline{*}}
\def \trir {\triangleright}
\def \tril {\triangleleft}
\def \Gdia {G_{{\rm diag}}}
\def \hs {\hspace{.2in}}
\def \Gdel {G_\delta}
\def \JJW {{}^{J^\prime}\!W^J}
\def \delJJW {{}^{\delta^{-1}(J^\prime)}\!W^{J}}
\def \Zw {Z_{\cc, \delta, w}}
\def \Zwp {Z_{\cc, \delta, w'}}
\def \Oxy {O_{\cc, x, y}}

\def \oG {\overline{G}}
\def \Zwm {Z_{\cc, \delta, w}^-}
\def \pe {\partial_{\emptyset}}
\def \WJmax {W^J_{\rm max}}
\def \WKmax {W^K_{\rm max}}
\def \WImax {W^I_{\rm max}}

\newtheorem{thm}{Theorem}[section]
\newtheorem{lem}{Lemma}[section]
\newtheorem{rmk}{Remark}[section]
\newtheorem{prop}{Proposition}[section]
\newtheorem{cor}{Corollary}[section]
\newtheorem{dfn}{Definition}[section]
\newtheorem{ex}{Example}[section]
\newtheorem{nota}{Notation}[section]
\newtheorem{nota-dfn}{Definition-Notation}[section]
\title[]{On intersections of certain partitions of a group compactification}
\author[Xuhua He]{Xuhua He}
\address{
Department of Mathematics \\
Hong Kong University of Science and Technology \\
Clear Water Bay\\
Hong Kong}
\email{maxhhe@ust.hk}
\author{Jiang-Hua Lu}
\address{
Department of Mathematics   \\
The University of Hong Kong \\
Pokfulam Road               \\
Hong Kong}
\email{jhlu@maths.hku.hk}
\date{}

\begin{abstract}
Let $G$ be  a connected semi-simple algebraic group of adjoint type over an algebraically closed field, and let
$\oG$ be the wonderful compactification of $G$.
For a fixed pair $(B, B^-)$ of opposite Borel subgroups of $G$, we look at intersections of Lusztig's $G$-stable pieces
and the $B^-\times B$-orbits in $\oG$, as well as intersections of $B \times B$-orbits and $B^- \times B^-$-orbits in
$\oG$. We give explicit conditions for such intersections to be non-empty, and in each case, we show that 
every non-empty intersection is smooth and irreducible, that the closure of the intersection is equal to 
the intersection of the closures, and that the non-empty intersections form a
strongly admissible partition of $\oG$.
\end{abstract}
\maketitle

\section{Introduction}
\subsection{}\label{subsec-partitions}
Let $Z$ be an irreducible algebraic variety over an algebraically closed field ${\bf k}$.
By a {\it partition} of $Z$ we mean a finite disjoint union $Z = \bigsqcup_{i \in \ci} X_i$ such that
each $X_i$ is a smooth irreducible locally closed subset of $Z$ and that
the closure 
of each $X_i$ in $Z$ is the union of some $X_{i'}$'s for $i' \in \ci$.

\subsection{}
Let $G$ be a connected semi-simple algebraic group of adjoint type
over an algebraically closed field ${\bf k}$. Regard $G$ as a 
$G \times G$ homogeneous space by the action
\[
(g_1, g_2) \cdot g = g_1 g g_2^{-1}, \hs g_1, g_2, g \in G.
\]
The De Concini-Procesi wonderful compactification $\oG$ of $G$ is a smooth $(G \times G)$-equivariant 
compactification of $G$ \cite{DP, DS}. 

Let $B$ and $B^-$ be a pair of opposite Borel subgroups of $G$. The partition of $\oG$ into the $B \times B$-orbits was studied in \cite{B} and \cite{Sp}. In his study of parabolic character sheaves on $\oG$ in \cite{L1, L2}, G. Lusztig introduced a decomposition of $\oG$ into finitely many $G$-stable pieces, where $G$ is identified with the diagonal $G_{\rm diag}$ of $G \times G$. It was later proved in \cite{H2} that Lusztig's $G$-stable pieces form a partition of $\oG$. 

This paper concerns with 

1) intersections of $B \times B$-orbits and $B^- \times B^-$-orbits in $\oG$,

2) intersections of the $G$-stable pieces
with $B^- \times B$-orbits in $\oG$.

Our motivation partially comes from Poisson geometry. Let $H = B \cap B^-$.
When ${\bf k}  = {\mathbb C}$,  
 there is \cite{e-l:cplx, LY2}
a natural $H \times H$-invariant Poisson structure $\Pi_1$ on $\oG$ whose $H \times H$-orbits of
symplectic leaves are the non-empty intersections of $B \times B$ and $B^- \times B^-$-orbits.
Similarly, there is natural 
$H_{\rm diag}$-invariant Poisson structure $\Pi_2$ on $\oG$ whose $H_{\rm diag}$-orbits of symplectic leaves are the non-empty 
intersections of $G_{\rm diag}$-orbits and $B^- \times B$-orbits. The restrictions
of $\Pi_1$ and $\Pi_2$ to $G \subset \oG$ are closely related to the
quantized universal enveloping algebra of the Lie algebra of $G$ and its dual (as a Hopf algebra).
See
\cite{e-l:grothendieck, KZ}. 

The closures of such intersections also appear in the study of
algebro-geometric properties of $\oG$. In the joint work \cite{HT2} of He and Thomsen, it was proved that in positive characteristics, there exists a Frobenius splitting on $\oG$ which compatibly splits all the nonempty intersections of the closures of $B \times B$-orbits and $B^- \times B^-$-orbits in $\oG$. In particular, all such closures are weakly normal and reduced. Moreover, the closure of a $B \times B$-orbit is globally F-regular in positive characteristic and is normal and Cohen-Macaulay for arbitrary characteristic. 

Later, in the joint work \cite{HT3} of He and Thomsen, it was proved that in positive characteristics, there exists a Frobenius splitting on $\oG$ which compatibly splits all the nonempty intersections of the closures of $G$-stable pieces and $B^- \times B$-orbits in $\oG$. In particular, all such closures are weakly normal and reduced. However, the closure of a $G$-stable piece is not normal in general \cite[No. 11.2]{HT3}.

\subsection{}
To state our results more precisely, we introduce some notation.
Let $N_G(H)$ be the normalizer of $H$ in $G$, and let $W=N_G(H)/H$ be the Weyl group.
Let $\Gamma$ be the set of simple roots  determined by the pair $(H, B)$.
For $J \subset \Gamma$, let $W_J$ be the subgroup of $W$ generated by  $J$,  and let $W^J \subset W$ 
the set of minimal length representatives of
$W/W_J$ in $W$.    If $J^\prime$ is a another subset of $\Gamma$, and $x \in W$,
let $\min(W_{J'} x W_J)$ and $\max(W_{J'} x W_J)$ be respectively
the unique minimal and maximal length elements in the double coset 
$W_{J'}x W_J$.

For $x, y \in W$, let $x \ast y \in W$ be such that $B (x \ast y) B$ is the unique
dense $(B, B)$-double coset in $BxByB$. 
The operation $\ast$ makes $W$ into a monoid  which will be denoted by $(W, \ast)$. See \cite{RS}.

Let $\delta \in {\rm Aut}(G)$ be  such that $\delta(H) = H$ and $\delta(B) = B$.
Let $G_\delta =\{(g, \delta(g)): \, g \in G\} \subset G \times G$ 
be the graph of $\delta$. We will in fact work with $G_\delta$-stable pieces in $\oG$
(see definition below). 

Recall that the $G \times G$-orbits in $\oG$ are in one to one correspondence with
subsets of $\Gamma$. For $J \subset \G$, let $Z_J$ be the corresponding $G \times G$-orbit in $\oG$. 
One has  
$\overline{Z_J} = \bigsqcup_{K \subset J} Z_K$,  and $\overline{Z_J}$ is smooth. See \cite{DP, DS}.
Let $h_J$ be a distinguished point in $Z_J$ (see $\S$\ref{subsec-recall-oG}).
For $w \in W^J$ and $(x, y) \in W^J \times W$, let 
\begin{align*}
&Z_{J, \delta, w}=\Gdel (B \times B)(w, 1)_\cdot h_J, \\
&[J, x, y]=(B \times B) (x, y)_\cdot h_J,\\
&[J, x, y]^{-, +}=(B^- \times B) (x, y)_\cdot h_J,\\
&[J, x, y]^{-, -}=(B^- \times B^-) (x, y)_\cdot h_J.
\end{align*}
The $Z_{J, \delta, w}$'s are called the $\Gdel$-stable pieces in $\oG$.
By \cite{H2, Sp}, one has the following partitions of $\oG$:
\begin{align}\label{eq-pppp}
\oG &=\bigsqcup_{J \subset \G, (x, y) \in W^J \times W} [J, x, y]=
\bigsqcup_{J \subset \G, (x, y) \in W^J \times W} [J, x, y]^{-, -} \\ 
\nonumber &=\bigsqcup_{J \subset \G, (x, y) \in W^J \times W} [J, x, y]^{-, +}
=\bigsqcup_{J \subset \G, w \in W^J} Z_{J, \d, w}.
\end{align}
For a subset $X$ of $\oG$, let 
$\overline{X}$ be the Zariski closure of $X$ in $\oG$.

We prove (see Proposition \ref{prop-Z-O-inter-oG}, Theorem \ref{thm-BB0}, and Theorem \ref{thm-inter-oG-2}) that for any $J \subset \Gamma$,
$w \in W^J$, and $(x, y), (u, v) \in W^J \times W$,

1) $[J, x, y] \cap [J, u, v]^{-, -} \neq \emptyset$ if and only if 
$x \le u$ and $v \le \max(yW_J),$ and in this case, $[J, x, y] \cap [J, u, v]^{-, -}$ is smooth and irreducible, and
\[
\overline{[J, x, y] \cap [J, u, v]^{-, -}} = \overline{[J, x, y]} \cap \overline{[J, u, v]^{-, -}}.
\]

2) $Z_{J, \delta, w} \cap [J, x, y]^{-, +} \neq \emptyset$ if and only if 
$\min (W_{J} \, \delta(w)) \le y^{-1} \ast \delta(x)$, and in this case, $Z_{J, \delta, w} \cap [J, x, y]^{-, +}$ is smooth and irreducible, and
\[
\overline{Z_{J, \delta, w} \cap [J, x, y]^{-, +}}=\overline{Z_{J, \delta, w}} \cap \overline{[J, x, y]^{-, +}}.
\]
Let
\begin{align*}
\cj &=\{(J, x, y, u, v): \;\;J \subset \Gamma, \,  (x, y), (u, v) \in W^J \times W,\\
& \hspace{1.3in} x \le u, \, v \le \max(yW_J)\},\\
\ck &= \{(J, w, x, y): \;\; J \subset \Gamma, \; (w, x, y) \in W^J \times W^J \times W, \\
& \hspace{1.3in} \min (W_{ J} \, \delta(w)) \le  y^{-1} \ast \delta(x)\}.
\end{align*}
One then has two more partitions of $\oG$:
\begin{align}\label{eq-pp}
\oG &= \bigsqcup_{(J, x, y, u, v) \in \cj} 
[J, x, y] \cap [J, u, v]^{-, -}\\
\nonumber &= \bigsqcup_{(J, w, x, y) \in \ck} Z_{J, \delta, w} \cap [J, x, y]^{-, +}
\end{align}

We introduce the notion of {\it admissible partitions} and {\it strongly admissible} partitions 
of $\oG$ (see Definition \ref{dfn-admissible})
and show that the six partitions in (\ref{eq-pppp}) and (\ref{eq-pp}) are all
strongly admissible (Proposition \ref{pr-Z-O-admissible}, Theorem \ref{thm-BB0}, 
and Theorem \ref{thm-inter-oG-2}). Moreover, the first two partitions in (\ref{eq-pppp}), as
well as the last two in (\ref{eq-pppp}), are shown to be {\it compatible}.
As consequences, we prove

1) if $J \subset \G$ and if
$X$ is  a subvariety of $Z_J$ appearing in any of the six partitions in (\ref{eq-pppp}) and (\ref{eq-pp}),
then for any $K \subset J$, $\overline{X} \cap Z_K \neq \emptyset$, and 
$\overline{X}$ and $\overline{Z_K}$ intersect properly in $\overline{Z_J}$. Moreover, we
describe the irreducible components of 
$\overline{X} \cap \overline{Z_K}$ in each case (Corollaries \ref{co-properly} and \ref{co-ZO-ZK}). 
This result for 
$\oG = \bigsqcup_{J \subset \G, (x, y) \in W^J \times W} [J, x, y]^{-, +}$ was also obtained by M. Brion \cite{B};

2) if $X = [J, x, y]$ and $Y = [K, u, v]^{-, -}$ with
$J, K \subset \G$, $(x, y) \in W^J \times W$, and $(u, v) 
\in W^K \times W$,  or if 
$X = Z_{J, \delta, w}$ and $Y = [K, x, y]^{-, +}$ with $J, K \subset \G, w\in W^J$ and $
(x, y) \in W^K \times W$, and if 
$\overline{X} \cap \overline{Y} \neq \emptyset$,  we show 
that $\overline{X}$ and $\overline{Y}$ intersect properly in $\overline{Z_{J \cup K}}$ (Corollary \ref{cor-inter-ZJK}).

In positive characteristic, let $G_F=\{(g, F(g)) : g \in G\}$ be the graph of Frobenius morphism $F: G \to G$. In $\S$\ref{sec-DL}, we study the intersection of the $G$-stable pieces with $G_F$-orbits. Such intersections include as a special case the Deligne-Lusztig varieties. 

Our discussions in this paper, and especially that in 
$\S$\ref{subsec-inter-general} and $\S$\ref{subsec-admissible}, also apply to intersections of 
$R$-stable pieces and $B \times B$-orbits, where $R$ is a certain class of connected subgroups of $G \times G$ as in \cite{LY},
as long as $R \cap (B \times B)$ is connected and that
${\rm Lie}(R) + {\rm Lie} (B \times B) ={\rm  Lie}(G \times G)$. 

\subsection{}\label{sec-notation}
We set up more notation for the rest of the paper. 

For $\alpha \in \Gamma$, let $U^\alpha$ be the one dimensional unipotent subgroup of $G$ defined by 
$\alpha$.
For a subset $J$ of $\Gamma$, let $P_J$ and $P_J^-$ be respectively the standard parabolic subgroups of $G$ 
determined by $J$  that contain $B$ and $B^-$, and let
$U_J$ and $U_J^-$ be respectively the 
uniradicals of $P_J$ and $P_J^-$.
Let $M_J=P_J \cap P_J^-$
be the common Levi factor of $P_J$ and $P_J^-$, and  let ${\rm Cen}(M_J)$ be the center of $M_J$.

The longest element in $W$ is denoted by $w_0$. If $J \subset \G$,
denote by $w_0^J$ the longest element in $W_J$, and  
let ${}^J\!W = \{w^{-1}: w \in W^J\}$.  If $J^\prime$ is a another subset of $\Gamma$,
let ${}^{J^\prime}\!W^J = {}^{J^\prime}\!W \cap W^J$.  

Throughout the paper, $\bigsqcup$ always means disjoint union.

\section{Intersections in $Z_{\cc}=(G \times G)/R_{\cc}$}\label{sec-inter}

\subsection{}\label{sec-2.1} Following \cite{LY}, an 
{\it admissible quadruple} for $G$  is a quadruple $\cc = (J, J^\prime, c, L)$, where
$J$ and $J^\prime$ are subsets of $\Gamma$, $c: J \to J^\prime$ is a bijective map preserving
the inner products between the simple roots, and $L$ is a connected closed subgroup of $M_{J} \times M_{J^\prime}$ of the form 
\[
L = \{(m, \, m^\prime) \in M_J \times M_{J^\prime}: \; \theta_c(m C) = m^\prime C^\prime\},
\]
with $C \subset {\rm Cen}(M_J)$ and $C^\prime \subset  {\rm Cen}(M_{J'})$ being closed subgroups and $\theta_c: 
M_J/C \to M_{J^\prime}/C^\prime$ a group isomorphism mapping $H/{\rm Cen}(M_J)$ to $H/{\rm Cen}(M_{J'})$ and
$U^\alpha$ to $U^{c(\alpha)}$ for every $\alpha \in J$.
For an admissible  quadruple $\cc = (J, J^\prime, c, L)$,
let
\begin{align}\label{eq-RC}
R_{\cc} &= L (U_J \times U_{J^\prime}) \subset P_J \times P_{J^\prime}. 
\end{align}
For example, $R_\cc = B \times B$ for $\cc = (\emptyset, \emptyset, {\rm Id}, H \times H)$ and
$R_\cc = \Gdia$ for $\cc=(\Gamma, \Gamma, {\rm Id}, \Gdia)$. For an admissible  quadruple $\cc = (J, J^\prime, c, L)$, let  
\[
Z_{\cc}=(G \times G)/R_{\cc}.
\]
When $G$ is of adjoint type, the $G \times G$-orbits in the De Concini-Procesi compactification $\oG$ of $G$ 
are all of the form $Z_\cc$ for some admissible quadruples $\cc$ (see $\S$\ref{subsec-recall-oG}).

\subsection{}\label{sec-2.2} For $(x, y) \in W^J \times W$, let
\begin{align*}
&[\cc, x, y]= (B \times B)(x, y)_\cdot R_{\cc} \subset Z_\cc,\\
&[\cc, x, y]^{-, +}= (B^- \times B)(x, y)_\cdot R_{\cc} \subset Z_\cc, \\
&[\cc, x, y]^{-, -}=(B^- \times B^-)(x, y)_\cdot R_\cc \subset Z_\cc.
\end{align*}
It follows from \cite{Sp} that 
\[
Z_\cc = \!\!\!\bigsqcup_{(x, y) \in W^J \times W}\!\!\! [\cc, x, y]
=\!\!\!\bigsqcup_{(x, y) \in W^J \times W} \!\!\![\cc, x, y]^{-, +}
=\!\!\!\bigsqcup_{(x, y) \in W^J \times W} \!\!\![\cc, x, y]^{-, -}
\]
are  the partitions of $Z_\cc$ by the $B \times B$, $B^-\times B$, and $B^- \times B^-$-orbits,
respectively.

\subsection{}\label{sec-2.3} Let $\delta$ be an automorphism of $G$ preserving both $H$ and $B$, and let 
\[
G_\delta =\{(g, \, \delta(g)): \, g \in G\} \subset G \times G
\]
be the graph of $\delta$. For  $w \in W^J$, let
\[
\Zw = \Gdel (B \times B)(w, 1)_\cdot R_{\cc}  \subset Z_\cc.
\]
The sets $\Zw$ for $w \in W^J$ will be called the {\it $\Gdel$-stable pieces} in $Z_\cc$. By \cite{H2, LY, Sp2},
each $\Zw$ is a locally closed smooth irreducible subset of $Z_\cc$, and 
\[
Z_\cc = \bigsqcup_{w \in W^J} \Zw
\]
is the partition of $Z_\cc$ by the $G_\d$-stable pieces.

\subsection{}\label{cl} We now recall the closure relations of the $B \times B$-orbits and $G_\d$-stable pieces in $Z_\cc$. 
For  $X \subset Z_\cc$, let $\overline{X}$ be the closure of $X$ in $Z_\cc$.

1) For $(x, y) \in W^J \times W$, $\overline{[\cc, x, y]}=\bigsqcup [\cc, x', y']$, where $(x', y')$ runs over elements in $W^J \times W$ such that $x' u \le x$ and $y' c(u) \le y$ for some $u \in W_J$. See \cite[Corollary 4.1]{LY}. 

2) For $w \in W^J$, $\overline{\Zw}=\bigsqcup Z_{\cc, \d, w'}$, where $w'$ runs over elements in $W^J$ such that $\d \i(c(u)) w' u \i \le w$ for some $u \in W_J$. See \cite[Corollary 5.9]{H3}. 

Using that facts that
\begin{align*}
[\cc, x, y]^{-, +}&=(w_0, 1) [\cc, w_0 x w_0^J, y w_0^{J'}], \\
[\cc, x, y]^{-, -}&=(w_0, w_0) [\cc, w_0 x w_0^J, w_0 y w_0^{J'}],
\end{align*}
one has the following variations of 1).

3) For $(x, y) \in W^J \times W$, $\overline{[\cc, x, y]^{-, +}}=\bigsqcup [\cc, x', y']^{-, +}$, where $(x', y')$ runs over elements in $W^J \times W$ such that $x' u \ge x w_0^J$ and $y' c(u) \le y w_0^{J'}$ for some $u \in W_J$.

4) For $(x, y) \in W^J \times W$, $\overline{[\cc, x, y]^{-, -}}=\bigsqcup [\cc, x', y']^{-, -}$, where $(x', y')$ runs over elements in $W^J \times W$ such that $x' u \ge x w_0^J$ and $y' c(u) \ge y w_0^{J'}$ for some $u \in W_J$.

\subsection{}\label{ast} Recall that the monoid operation $\ast$ on $W$ is defined by $\overline{B (x \ast y) B}\!\! = \overline{BxByB}$
for $x, y \in W$.
Similarly, for $x, y \in W$, define 
$x \trir y \in W$ and $x \tril y \in W$ by
\[
\overline{BxByB^-} = \overline{B (x \trir y) B^-} \hs \mbox{and} \hs 
\overline{B^-xByB} = \overline{B^- (x \tril y) B}.
\]
Then 
\[
(W, \ast) \times W \longrightarrow  W: \; (x, y) \longmapsto x \trir y, \hs x, y \in W
\]
is a left monoid action of $(W, \ast)$ on $W$, and 
\[
W \times (W, \ast) \longrightarrow W: \; (x, y) \longmapsto x \tril y, \hs x, y \in W
\]
is a right monoidal action of $(W, \ast)$ on $W$. More properties of $\ast, \trir$ and $\tril$ are 
reviewed in the Appendix.

\subsection{}\label{subsec-inter-1}
We now determine when the intersection of a $B \times B$-orbit and a $B^- \times B^-$-orbit in 
$Z_\cc$ is non-empty. 

\begin{prop}\label{pr1-inter}
For any $(x, y), (u, v) \in W^J \times W$, the following conditions are equivalent:

1) $[\cc, x, y] \cap [\cc, u, v]^{-, -} \neq \emptyset$,

2) $u \le x$ and $\min(v W_{J'}) \le y$,

3) $u \le x$ and $v \le \max(y W_{J'})$.
\end{prop}

{\bf Proof.} 
Using the facts that $x, u \in W^J$, it is easy to see that 
\begin{align*}
&[\cc, x, y]=(B \times B) (x, y) (B \times B)_\cdot R_\cc,\\
&[\cc, u, v]^{-, -}=(B^- \times B^-) (uw_0^J,vw_0^{J'}) (B \times B)_\cdot R_\cc.
\end{align*}
Thus $[\cc, x, y] \cap [\cc, u, v]^{-, -} \neq \emptyset$ if and only if 
$$(B x B, B y B) \cap \bigl( (B^- \times B^-) (uw_0^J,vw_0^{J'}) (B \times B) R_\cc (B \times B) \bigr) \neq \emptyset.$$
Since
\begin{equation}\label{eq-BRB}
(B \times B) R_{\cc} (B \times B) = \bigcup_{z \in W_J} (B \times B) (z, c(z)) (B \times B),
\end{equation}
$[\cc, x, y] \cap [\cc, u, v]^{-, -} \neq \emptyset$ if and only if 
\begin{equation}\label{eq-BBBB-0}
(BxB, \;B y B) \cap (B^- uw_0^J BzB, \; B^- v w_0^{J'}B c(z) B)  \neq \emptyset
\end{equation}
for some $z \in W_J$. By Lemma \ref{lem-facts3} and Lemma \ref{lem-facts1} in the Appendix, 
(\ref{eq-BBBB-0}) is the same as
\[
u w_0^J\le x \ast z^{-1} \hs \mbox{and} \hs v w_0^{J'}\le y \ast c(z)^{-1} \hs\
\mbox{for some} \; z \in W_J.
\]
Since for any $z \in W_J$, $x \ast z \le \max(xW_J)$ and $ y \ast c(z) \le \max(y W_{J'})$
and both inequalities become equalities when $z = w_0^J$, (\ref{eq-BBBB-0}) is equivalent to
$uw_0^J \le \max(xW_J)$ and $ vw_0^{J'} \le \max(y W_{J'})$ which, by Lemma \ref{lem-x} in the Appendix, 
are in turn equivalent to
$u \le x $ and $\min(v W_{J'}) \le y$, or $u \le x$ and $v \le \max(y W_{J'})$.  \qed

\begin{ex}
{\em When $R_\cc = G_{\rm diag}$ and $Z_\cc$ is identified with $G$, 
the intersections in Proposition \ref{pr1-inter} are of the form $ByB \cap B^- w B^-$
for $y, w \in W$, and are called {\it double Bruhat cells} \cite{fz-double-positive}. It is well-known 
(see, for example, \cite{fz-double-positive}) that
the intersection
$(ByB) \cap (B^- w B^-)$ is non-empty for all $y, w \in W$, which can also be seen from Proposition \ref{pr1-inter}.}
\end{ex}

\

\subsection{}\label{subsec-inter-2}
We now determine when the intersection of a $\Gdel$-stable piece and a $B^- \times B$-orbit in 
$Z_\cc$ is nonempty. 

\begin{prop}\label{pr-inter}
For  $w\in W^J$ and $(x, y) \in W^J \times W$, the following conditions are equivalent: 

1) $\Zw \cap [\cc, x, y]^{-, +} \neq \emptyset$, 

2) $y^{-1} \trir \delta(x) \le \max(W_{J'} \, \delta(w))$,

3) $\min (W_{J'} (y^{-1} \trir \delta(x))) \le \delta(w)$.
\end{prop}

\noindent
{\bf Proof.} Using the facts that $w, x \in W^J$, it is easy to see that
\begin{align*}
&\Zw = \Gdel (B \times B)(w, 1)(B \times B)_\cdot R_{\cc},\\
&[\cc, x, y]^{-, +}  = (B^- \times B)(xw_0^J, yw_0^{J'})(B \times B)_\cdot R_{\cc}.
\end{align*}
Thus $\Zw \cap [\cc, x, y]^{-, +} \neq \emptyset$ if and only if
\[
\Gdel \cap \left(
(B^- \times B)(xw_0^J, yw_0^{J'})(B \times B)R_{\cc}(B \times B)(w \i, 1)(B \times B)
\right) \neq \emptyset.
\]
By (\ref{eq-BRB}), $\Zw \cap [\cc, x, y]^{-, +} \neq \emptyset$ if and only if 
\[
\Gdel \cap \left(
(B^- \times B)(xw_0^J, yw_0^{J'})(B \times B)(z, c(z))(B \times B)(w \i, 1)(B \times B)
\right) \neq \emptyset
\]
for some $z \in W_J$, which is equivalent to
\begin{equation}\label{eq-1}
(B^- \delta(xw_0^J) B \delta(z) B \delta(w \i) B) \cap (B yw_0^{J'} B c(z) B) \neq \emptyset.
\end{equation}
Since $l(z w^{-1}) = l(z) + l(w^{-1})$ for every $z \in W_J$,  (\ref{eq-1}) is equivalent to
\[
(B^- \delta(xw_0^J) B \delta(z w \i) B) \cap (B yw_0^{J'} B c(z) B) \neq \emptyset
\]
which, by Lemma  \ref{lem-facts3} and Lemma \ref{lem-facts1} in the Appendix, is equivalent to
\begin{equation}\label{eq-2}
\delta(x w_0^J) \le (y w_0^{J'}) \ast  c(z) \ast \delta(w) \ast \delta(z^{-1})
\end{equation}
for some $z \in W_J$. Since by Lemma \ref{lem-tri-min} in the Appendix,
\begin{align*}
(y w_0^{J'}) \ast  c(z) \ast \delta(w) \ast \delta(z^{-1})
& \le (y w_0^{J'}) \ast c(w_0^J) \ast \delta(w) \ast \delta(w_0^J)\\
& = y \ast \max(W_{J'} \, \delta(w) \, W_{\delta(J)})
\end{align*}
for any $z \in W_J$ with equality holds  when $z=w_0^J$,
(\ref{eq-2}) is equivalent to 
\begin{equation}\label{eq-3}
\delta(x w_0^J) \le  y \ast \max(W_{J'} \, \delta(w) \, W_{\delta(J)}).
\end{equation}
Clearly (\ref{eq-3}) leads to 
\begin{equation}\label{eq-4}
\delta(x) \le y \ast \max(W_{J'} \, \delta(w) \, W_{\delta(J)}).
\end{equation}
Conversely, if (\ref{eq-4}) holds, then 
\begin{align*}
\delta(x w_0^J)& = \delta(x) \ast w_{0}^{\delta(J)} \le y \ast \max(W_{J'} \, \delta(w) \, W_{\delta(J)}) \ast w_0^{\delta(J)} \\
& = y \ast \max(W_{J'} \, \delta(w) \, W_{\delta(J)}).
\end{align*}
Thus (\ref{eq-3}) is equivalent to (\ref{eq-4}), which, by 
Lemma \ref{lem-facts1} in the Appendix, is equivalent to
\begin{equation}\label{eq-5}
y^{-1} \trir \delta(x) \le \max(W_{J'} \, \delta(w) \, W_{\delta(J)}).
\end{equation}
Since $y^{-1} \trir \delta(x) \in W^{\delta(J)}$ by Lemma \ref{lem-tri-min} in the Appendix, and since 
\[
\max(W_{J'} \, \delta(w) \, W_{\delta(J)}) = \max(W_{J'} \, \delta(w))\ast w_0^{\delta(J)},
\] 
(\ref{eq-5}) is equivalent to $y^{-1} \trir \delta(x) \le \max(W_{J'} \, \delta(w))$
by Lemma \ref{lem-x} in the Appendix. The equivalence of 2) and 3) also follows from Lemma \ref{lem-x} in the Appendix. \qed

\

\subsection{}\label{subsec-consequences}
We now discuss some consequences of the results in $\S$\ref{subsec-inter-1} and $\S$\ref{subsec-inter-2}.

\begin{cor}
Let $J \subset \G$. For $(x, y) \in W^J \times W$. Set $$w_{x, y}=\min(W_{\d \i(J')} (\d \i(y \i) \trir x)) \in {}^{\d \i(J')} W^J.$$ 
Then for $w \in W^J$, 
\[
\Zw \cap [\cc, x, y]^{-, +} \neq \emptyset\hs \mbox{iff} \hs Z_{\cc, \d, w_{x, y}} \subset \overline{\Zw}.
\] 
\end{cor}

{\bf Proof.} If $\Zw \cap [\cc, x, y]^{-, +} \neq \emptyset$, by Proposition \ref{pr-inter}, $w_{x, y} \le w$. By $\S$\ref{cl}, 2), $Z_{\cc, \d, w_{x, y}} \subset \overline{\Zw}$. On the other hand, if $Z_{\cc, \d, w_{x, y}} \subset \overline{\Zw}$, then there exists $u \in W_J$ such that $\d \i(c(u)) w_{x, y} u \i \le w$. Since $w_{x, y} \in {}^{\d \i(J')} W^J$, $w_{x, y} \le \d \i(c(u)) w_{x, y} u \i \le w$. By Proposition \ref{pr-inter}, $\Zw \cap [\cc, x, y]^{-, +} \neq \emptyset$. \qed

\begin{prop}\label{prop-pi}
Let $\pi: Z_\cc \to  (G \times G)/(P_J \times P_{J'})$ be the natural projection induced by the inclusion $R_\cc \subset P_J \times P_{J'}$. Then for any $w \in W^J$ and $(x, y) \in W^J \times W$, \[
\Zw \cap [\cc, x, y]^{-, +} \neq \emptyset \; \hs \mbox{iff} \hs \;
\pi(\Zw) \cap \pi([\cc, x, y]^{-, +}) \neq \emptyset.
\]
\end{prop}

{\bf Proof.}
Clearly $\Zw \cap [\cc, x, y]^{-, +} \neq \emptyset$ implies that
$\pi(\Zw) \cap \pi([\cc, x, y]^{-, +}) \neq \emptyset.$ Assume now that $\pi(\Zw) \cap \pi([\cc, x, y]^{-, +}) \neq \emptyset.$ Let $y'=\min(y W_{J'}) \in W^{J'}$ and $w'=\min(W_{\d \i(J')} w) \in {}^{\d \i(J')} W^J$. Then 
\begin{align*} 
&\pi(\Zw)=\Gdel (w', 1)(P_J \times P_{J^\prime}), \\ 
&\pi([\cc, x, y]^{-, +})=(B^- \times B)(x, y') (P_J \times P_{J^\prime}).
\end{align*} 
By definition, $\max(W_{J'} \d(w))=\max(W_{J'} \d(w'))$. By Lemma \ref{2} in the Appendix, $y \i \trir \d(x) \le (y') \i \trir \d(x)$. Now $\Zw \cap [\cc, x, y]^{-, +} \neq \emptyset$ follows from Proposition \ref{pr-inter} and the following Lemma \ref{lem-PJJ}.  \qed

\begin{lem}\label{lem-PJJ}
For $w \in \delJJW$ and $(x, y)  \in W^J \times W^{J^\prime}$,
\[
\Gdel (w, 1)(P_J \times P_{J^\prime}) \cap (B^- \times B)(x, y) (P_J \times P_{J^\prime})
\neq \emptyset
\]
if and only if  $y^{-1} \trir \delta(x) \le \max(W_{J'} \, \delta(w)).$
\end{lem}

{\bf Proof.}
First note that 
\[
\Gdel (w, 1) (P_J \times P_{J^\prime}) = \Gdel (B \times B) (w, 1) (P_J \times P_{J^\prime}).
\]
Thus $\Gdel (w, 1) (P_J \times P_{J^\prime}) \cap (B^- \times B)(x, y)(P_J \times P_{J^\prime})
\neq \emptyset$ if and only if
\begin{equation}\label{eq-Gdel-JJ}
\Gdel \cap (B^- \times B)(x, y) (P_J \times P_{J^\prime}) (w^{-1}, 1) (B \times B) \neq \emptyset.
\end{equation}
Using  $P_J \times P_{J'} = \cup_{z \in W_J, z' \in W_{J'}} (B \times B)(z, z') (B \times B)$ and the fact that
$BzBw^{-1}B = Bzw^{-1}B$ for any $z \in W_J$, one sees that
(\ref{eq-Gdel-JJ}) is equivalent to 
\[
\bigcup_{z \in W_J, z' \in W_{J'}}  \Gdel \cap (B^- x B zw^{-1} B \times B y B z'B) \neq \emptyset,
\]
or $(B^- \delta(x) B \delta(zw^{-1})B) \cap (B y B z' B) \neq \emptyset$ for some $z \in W_J$ and $z' \in W_{J'}$,
which, by Lemma \ref{lem-facts3} and Lemma \ref{lem-facts1} in the Appendix, is in turn equivalent to 
\[
y^{-1} \trir \delta(x) \le z' \ast \delta(w) \ast \delta(z^{-1})
\hs \mbox{for some} \; z \in W_J, \, z' \in W_{J'}.
\]
Since for any $z \in W_J$ and $z' \in W_{J'}$,
\[
z' \ast \delta(w) \ast \delta(z^{-1}) \le w_0^{ J'} \ast \delta(w) \ast w_0^{\delta(J)}
=\max(W_{J'} \, \delta(w) \, W_{\delta(J)})
\]
with  equality holds when $z=w_0^J$ and $z' = w_0^{J'}$, 
(\ref{eq-Gdel-JJ}) is equivalent to 
\begin{equation}\label{eq-6}
y^{-1} \trir \delta(x) \le \max(W_{J'} \, \delta(w) \, W_{\delta(J)}).
\end{equation}
Since $\delta(x) \in W^{\delta(J)}$, it follows from Lemma \ref{lem-tri-min} 
in the Appendix that (\ref{eq-6}) is
equivalent to $y^{-1} \trir \delta(x) \le \max(W_{J'} \, \delta(w))$. \qed


\subsection{}\label{subsec-inter-general} To study the geometry and closures of the non-empty intersections
in $\S$\ref{subsec-inter-1} and $\S$\ref{subsec-inter-2}, 
we first recall some elementary facts on intersections of subvarieties in 
an algebraic variety.

The following Lemma \ref{lem-Richardson}  is
 a generalization of 
\cite[Corollary 1.5]{R} of Richardson. Our proof of Lemma \ref{lem-Richardson} is essentially 
the same as that of \cite[Theorem 1.4]{R}.

\begin{lem}\label{lem-Richardson}
Let $A$ be a connected algebraic group and let $H, K$ and $L$ be closed connected subgroups of $A$. Assume that
$H \cap K$ is connected and that ${\rm Lie}(H) + {\rm Lie}(K) = {\rm Lie}(A)$. Let $Y$ be an 
irreducible subvariety of $A/L$ such that $HY \subset A/L$ is smooth. Then for any $K$-orbit  $O$ in $A/L$ 
such that $(HY) \cap O \neq \emptyset$, $HY$ and $O$ intersect transversally in $A/L$ and
$HY \cap O$ is a smooth irreducible subvariety of $A/L$ with
\[
\dim ((HY) \cap O) = \dim HY + \dim O - \dim A/L.
\]
\end{lem}

{\bf Proof.}  Since $HY$ is a union of $H$-orbits in $A/L$, it follows from
\cite[Corollary 1.5]{R} and \cite[Proposition 1.2]{R}
that $HY$ and $O$ intersect transversally and that the 
intersection $(HY) \cap O$ is smooth. Moreover, each irreducible component of $(HY) \cap O$ has 
dimensional equal to $\dim HY + \dim O - \dim A/L$.

It remains to show that $(HY) \cap O$ is irreducible.
Fix an $x \in O$ and consider the diagram 
\[
O \stackrel{p}{\longleftarrow} H \times K \stackrel{m}{\longrightarrow} A 
\stackrel{q}{\longrightarrow} A/L,
\] 
where $p(h, k)=k x$, $m(h, k)=h \i k$, and $q(a) = ax$ for $h \in H, k \in K$, and $a \in A$.
Let
\[
E= \{(h, k) \in H \times K: \; h^{-1} k x \in Y\} \subset H \times K.
\]
Then $(HY) \cap O=p(E)$, so it is enough to show that $E$ is irreducible.

Since $L$ is connected and $Y \subset A/L$ is irreducible, $q^{-1}(Y) \subset A$ is irreducible by \cite[Lemma 1.3]{R}. As in the proof of \cite[Theorem 1.4]{R}, $HK$ is open in $A$, so
$HK \cap q \i(Y)$ is an irreducible subvariety of $HK$. 
The map $m$ induces an isomorphism $m: (H\times K)/D \rightarrow HK$, where $D = \{(g, g): g \in H \cap K\}$. Let
$\nu: H \times K \to (H \times K)/D$ be the natural projection. Since $D$ is connected, by \cite[Lemma 1.3]{R}, 
$E=\nu^{-1}(m \i(HK \cap q \i(Y)))$ is also irreducible. \qed

\

The following Lemma \ref{le-general} is useful in determining the irreducible components of
intersections of algebraic varieties and will be used several times in the paper. 

\begin{lem}\label{le-general}
Let $Y$ be an algebraic variety over an algebraically closed field ${\bf k}$. Suppose that $l \geq 0$ is an integer such that
every irreducible component of $Y$ has dimension at least $l$, and suppose that 
$Y = \bigsqcup_{k \in \ck} Y_k$  is a finite disjoint union,
where each $Y_k$ is an irreducible subvariety of $Y$ with $\dim Y_k \le l$.
Then the irreducible components of $Y$ are precisely the closures $\overline{Y_k}$ of those $Y_k$'s, where  $k \in \ck$ and $\dim Y_k = l$. 
\end{lem}

{\bf Proof.} Let $S$ be any irreducible component of $Y$. Then 
\[
S = \bigcup_{k \in \ck: S \cap Y_k \neq \emptyset} \overline{S \cap Y_k}.
\]
Since $S$ is irreducible, $S = \overline{S \cap Y_k} \subset \overline{Y_k}$ for some $k \in \ck$. Since $\overline{Y_k}$
is irreducible, $S = \overline{Y_k}$, and it follows from the dimension assumptions that $\dim Y_k = l$. Since the $Y_k$'s 
are pair-wise disjoint, the closures $\overline{Y_k}$ with $\dim Y_k = l$ are pair-wise distinct irreducible components. \qed

\begin{lem}\label{lem-closure-inter}
Let $Z$ be a smooth irreducible algebraic variety and let 
\[
Z = \bigsqcup_{i \in \ci} X_i = \bigsqcup_{j \in \cj} Y_j
\]
be two partitions of $Z$ such that each non-empty intersection $X_i \cap Y_j$ is transversal and irreducible. Then
for any $(i, j) \in \ci \times \cj$, $X_i \cap Y_j \neq \emptyset$ if and only if
$\overline{X_i} \cap \overline{Y_j} \neq \emptyset$, and in this case,
\[
\overline{X_i \cap Y_j} = \overline{X_i} \cap \overline{Y_j}.
\]
In particular, $Z = \bigsqcup_{(i, j) \in \ck} (X_i \cap Y_j)$
is again a partition of $Z$. Here $\ck=\{(i, j) \in \ci \times \cj: X_i \cap Y_j \neq \emptyset\}$. 
\end{lem}

{\bf Proof.} 
Let $(i, j)  \in \ci \times \cj$ be such that $\overline{X_i} \cap \overline{Y_j} \neq \emptyset$, and let
\[
\ck_{ij} = \{(i', j') \in \ci \times \cj: \,X_{i'} \subset \overline{X_i}, \,\, Y_{j'} \subset \overline{Y_j}, \, 
X_{i'} \cap Y_{j'} \neq \emptyset\}.
\]
Then
\[
\overline{X_i} \cap \overline{Y_j} = \bigsqcup_{(i', j') \in \ck_{ij}} X_{i'} \cap Y_{j'}
\]
is a disjoint union. By \cite[Page 222]{Hr}, every irreducible component of $\overline{X_i} \cap \overline{Y_j}$
has dimension at least $\dim X_i + \dim Y_j - \dim Z$. On the other hand, for any $(i', j') \in \ck_{ij}$,
\begin{align}\label{eq-dims}
\dim X_{i'} \cap Y_{j'} &= 
\dim X_{i'} + \dim Y_{j'} - \dim Z \\
\nonumber
 &\le \dim X_{i} + \dim Y_{j} - \dim Z.
\end{align}
Since $\overline{X_i}$ is irreducible, $\dim X_{i'} < \dim X_i$ for
any $i' \in \ci$ such that $X_{i'} \subset 
\overline{X_i}$ and $i' \neq i$. Similarly, $\dim Y_{j'} < \dim Y_j$ for
any $j' \in \cj$ such that $Y_{j'} \subset 
\overline{Y_j}$ and $j' \neq j$. 
Thus the inequality in (\ref{eq-dims}) is an equality if and only if $(i', j') = (i, j)$.
By Lemma \ref{le-general}, 
$X_i \cap Y_j \neq \emptyset$ and $\overline{X_i \cap Y_j} = \overline{X_i} \cap \overline{Y_j}$. \qed

\begin{thm}\label{thm-closure-inter}
Let $w \in W^J$ and $(x, y), (u, v) \in W^J \times W$. Then 

1) $\Zw \cap [\cc, x, y]^{-, +} \neq \emptyset$ if and only if 
$\overline{\Zw} \cap \overline{[\cc, x, y]^{-, +}}
\neq \emptyset$. In this case, $\Zw$ and $[\cc, x, y]^{-, +}$ intersects transversally in 
$Z_\cc$, the intersection is smooth and irreducible, and 
\[
\overline{\Zw \cap [\cc, x, y]^{-, +}} = \overline{\Zw} \cap \overline{[\cc, x, y]^{-, +}}.
\]

2) $[\cc, x, y] \cap [\cc, u, v]^{-, -} \neq \emptyset$ if and only if $\overline{[\cc, x, y]} \cap \overline{[\cc, u, v]^{-, -}}
\neq \emptyset$. In this case, $[\cc, x, y]$ and $[\cc, u, v]^{-, -}$ intersects 
transversally in $Z_\cc$, the intersection is smooth and irreducible, and 
\[
\overline{[\cc, x, y] \cap [\cc, u, v]^{-, -}} =\overline{[\cc, x, y]} \cap \overline{[\cc, u, v]^{-, -}}.
\]
\end{thm}

{\bf Proof.} Since $R_\cc$ is connected, $\Gdel \cap (B^- \times B)$ is connected, and 
\[
{\rm Lie}(\Gdel) + {\rm Lie}(B^- \times B) = {\rm Lie}(G \times G),
\] 
Lemma \ref{lem-Richardson}  applies. By taking $A = G \times G$, 
\[
H = G_{\delta}, \; \; K = B^- \times B, \; \; L = R_\cc, \;\; \mbox{and} \;\; 
Y = (B \times B)(w, 1)_\cdot R_\cc
\]
in Lemma \ref{lem-Richardson}, one sees that when $\Zw \cap [\cc, x, y]^{-, +} \neq \emptyset$, $\Zw$ 
and $[\cc, x, y]^{-, +}$ intersect transversally in $Z_\cc$, and that  the intersection 
$\Zw \cap [\cc, x, y]^{-, +}$ is smooth and irreducible.
By applying 
Lemma  \ref{lem-closure-inter} to the two
partitions
\[
Z_\cc = \bigsqcup_{w \in W^J} \Zw = \bigsqcup_{(x, y) \in W^J \times W} (B^- \times B)(x, y)_\cdot R_\cc
\]
of $Z_\cc$, one proves part 1). 
Part 2) can be proved in the same way.  \qed 

\begin{rmk}
{\em In both 1) and 2) in Theorem \ref{thm-closure-inter}, the fact that the intersection is non-empty
if and only if the intersection of the closures is non-empty can also been obtained using $\S$\ref{cl} 
and Proposition \ref{pr1-inter} and  Proposition \ref{pr-inter}. However, the proof we give is more conceptual.} 
\end{rmk}
\section{Intersections in $\oG$}

\subsection{}\label{subsec-recall-oG}
Let $G$ be a connected semi-simple adjoint group and $\oG$ be the De Concini-Procesi compactification. 
It is well-known that the $G \times G$-orbits in $\oG$ are in one to one correspondence with
subsets of $\Gamma$. For $J \subset \G$, let $Z_J$ be the corresponding $G \times G$-orbit in $\oG$. 
A distinguished point $h_J \in Z_J$ can be chosen such that the stabilizer subgroup of $G \times G$ at $h_J$ is
\[
R_J^- \stackrel{{\rm def}}{=} (U_J^- \times U_J)\{(m_1, m_2) \in M_J \times M_J: \; \pi_J(m_1) = \pi_J(m_2)\},
\]
where $\pi_J: M_J \to M_J/{\rm Cen}(M_J)$ is the natural projection and
${\rm Cen}(M_J)$ is the center of $M_J$.

For $J \subset \G$, let $\cc_J=(J^*, J, c, L)$, where 
$J^*=-w_0(J)$, $c=(w_0w_0^J)^{-1}$, and 
\[
L=\{(\dot{w}_0 \dot{w}_0^J m_1 (\dot{w}_0 \dot{w}_0^J) \i, m_2): m_1, m_2 \in M_J, \pi_J(m_1)=\pi_J(m_2)\}
\]
with $\dot{w}_0$ and $\dot{w}_0^J$ being any representatives of
$w_0$ and $w_0^J$ in $N_G(H)$.
By \cite[Section 5]{LY}, $\cc_J$ is an admissible quadruple for $G$, 
and 
\[
R_{\cc_J} = (\dot{w}_0 \dot{w}_0^J, 1) R^-_J (\dot{w}_0 \dot{w}_0^J, 1)^{-1}.
\]
One thus has the isomorphism
\begin{equation}\label{eq-ZZ}
Z_J \longrightarrow  Z_{\cc_J}: \; \; \; (g, g')_\cdot h_J \longmapsto 
(g w_0^J w_0, g')_\cdot R_{\cc_J}, \hs g, g^\prime \in G.
\end{equation}

\subsection{}\label{dim} For $J \subset \Gamma$ and $(x, y) \in W^J \times W$, let 
\begin{align}\label{J-0}
&[J, x, y]=(B \times B) (x, y)_\cdot h_J,\\
\label{J-1}&[J, x, y]^{-, +}=(B^- \times B) (x, y)_\cdot h_J=(w_0, 1) [J, \;w_0 x w_0^J,\; y w_0^J], \\
\label{J-2}&[J, x, y]^{-, -}=(B^- \times B^-) (x, y)_\cdot h_J=(w_0, w_0) [J,\; w_0 x w_0^J,\; w_0 y w_0^J].
\end{align}
For $J \subset \G$ and $w \in W^J$, let 
$$Z_{J, \delta, w}=\Gdel (B \times B)(w, 1)_\cdot h_J.$$
The $Z_{J, \delta, w}$'s will be called the $\Gdel$-stable pieces in $\oG$.
By \cite{H2, Sp}, one has the following partitions of $\oG$:
\begin{align*}
\oG &=\bigsqcup_{J \subset \G, (x, y) \in W^J \times W} [J, x, y]=
\bigsqcup_{J \subset \G, (x, y) \in W^J \times W} [J, x, y]^{-, -} \\ 
&=\bigsqcup_{J \subset \G, (x, y) \in W^J \times W} [J, x, y]^{-, +}
=\bigsqcup_{J \subset \G, w \in W^J} Z_{J, \d, w}.
\end{align*}

For a an irreducible subvariety $X \subset Z_J$, 
let ${\rm Codim}_{Z_J}(X)$ be the codimension of $X$ in $Z_J$. 
Let $l$ be the length function of $W$. One has, for any $J \subset \Gamma$
and $w \in W^J$,

1) $\dim Z_J = \dim G - \dim {\rm Cen}(M_J)=\dim G - |\G| + |J|$. See \cite{DP}. 

2) ${\rm Codim}_{Z_J}[J, x, y] = l(w_0)+l(x) -l(y)$. See \cite[Lemma 1.3]{Sp}.

3) ${\rm Codim}_{Z_J}Z_{J, \delta, w} = l(w)$. See \cite[Section 8]{L2}. 

\noindent
By (\ref{J-1}) and (\ref{J-2}),
one also has 

4) ${\rm Codim}_{Z_J}[J, x, y]^{-, +} =2 l(w_0)-l(x w_0^J)-l(y w_0^J)$.

5) ${\rm Codim}_{Z_J}[J, x, y]^{-, -} =l(w_0)-l(x w_0^J)+l(y w_0^J)$.

\subsection{}\label{cl'}

The closure of a $G \times G$-orbit is described in \cite{DP, DS} as follows. 

1) For $J \subset \Gamma$, $\overline{Z_J} = \bigsqcup_{K \subset J} Z_K$ is a smooth subvariety of $\oG$.

The closure of a $B \times B$-orbit is described in \cite[Proposition 2.4]{Sp}, and the following simplified version in 2) 
is found in \cite[Proposition 6.3]{HT2} and \cite[Example 1.3]{LY}.
The following 
3) and 4) are obtained using (\ref{J-1}) and (\ref{J-2}).

2) For $J \subset \G$ and $(x, y) \in W^J \times W$,  $\overline{[J, x, y]}=\bigsqcup [K, x', y']$, 
where $K \subset J$, $(x', y') \in W^{K} \times W$ and there exists $u \in W_J$ such that $x u \le x', y' \le y u$. 

3) For $J \subset \G$ and $(x, y) \in W^J \times W$, $\overline{[J, x, y]^{-, +}}=\bigsqcup [K, x', y']^{-, +}$, 
where $K \subset J$, $(x', y') \in W^{K} \times W$ and there exists $u \in W_J$ such that $x' w_0^{K} \le x u, y' w_0^{K} \le y u$. 

4) For $J \subset \G$ and $(x, y) \in W^J \times W$,   $\overline{[J, x, y]^{-, -}}=\bigsqcup [K, x', y']^{-, -}$, 
where $K \subset J$, $(x', y') \in W^{K} \times W$ and there exists $u \in W_J$ such that $x' w_0^{K} \le x u, y' w_0^{K} \ge y u$.

For $J \subset \Gamma$ and $w \in W$, let 
\[
C_J(w) = \{\delta^{-1}(u) wu^{-1}: u \in W_J\},
\]
and denote by ${\rm Min}(C_J(w))$ the set of minimal length elements in $C_J(w)$.
The closure of a $G_\d$-stable piece is described in \cite[Sections 3 and 4]{H2} as follows:

5) For $J \subset \Gamma$ and $w \in W^J$, $\overline{Z_{J, \delta, w}} =\bigsqcup Z_{K, \delta, w'}$, 
where $K \subset J, w' \in W^{K}$, and $w'\ge w_1$ for some $w_1 \in {\rm Min}(C_J(w))$.

\
\subsection{} We can now prove our first main result in this paper.

\begin{prop}\label{prop-Z-O-inter-oG}
For $J \subset \Gamma$, $w \in W^J$, and $(x, y), (u, v) \in W^J \times W$, 
\begin{align}\label{J-iff-1}
 &[J, x, y] \cap [J, u, v]^{-, -} \neq \emptyset \; \; \mbox{iff} \;\;\;
x \le u, \, v \le \max(yW_J), \\
\label{J-iff-2}
& Z_{J, \delta, w} \cap [J, x, y]^{-, +} \neq \emptyset \; \; \mbox{iff}\;\;\;
\min (W_{J} \, \delta(w)) \le y^{-1} \ast \delta(x).
\end{align}
\end{prop}

{\bf Proof.} 
Let $\cc_J$ be as in $\S$\ref{subsec-recall-oG} and recall the isomorphism 
$Z_J \to Z_{\cc_J}$ in (\ref{eq-ZZ}). Since $W^{J^*} = W^J w_0^J w_0$, one has
\[
[J, x, y] \cap [J, u, v]^{-, -} \neq \emptyset
\; \; \mbox{iff} \;\;\; [\cc_J, x w_0^J w_0, y] \cap [\cc_J, u w_0^J w_0, v]^{-, -} \neq \emptyset,
\]
which, by Proposition \ref{pr1-inter}, is equivalent to $u w_0^J w_0 \le x w_0^J w_0$ and $v \le \max(y W_J)$. 
Note that $u w_0^J w_0 \le x w_0^J w_0$ if and only if $u w_0^J \ge x w_0^J$, which is equivalent to $u \ge x$ since $x, u \in W^J$. Thus 
(\ref{J-iff-1}) is proved. 

Similarly, $Z_{J, \d, w} \cap [J, x, y]^{-, +} \neq \emptyset$ iff $Z_{\cc, \delta, w w_0^J w_0} 
\cap [\cc, x w_0^J w_0, y]^{-, +} \neq \emptyset$, which, by Proposition \ref{pr-inter}, is equivalent to 
\begin{equation}\label{eq-7}
y^{-1} \trir (\delta(x) w_0^{\d(J)} w_0) \le \max(W_J \delta(w) w_0^{\delta(J)} w_0).
\end{equation}
By Lemma \ref{lem-x} and Lemma \ref{lem-facts1} in the Appendix, (\ref{eq-7}) is equivalent to 
\[
\min (W_J \delta(w) w_0^{\delta(J)}) \le y^{-1} \ast (\delta(x) w_0^{\d(J)}),
\]
which is in turn equivalent to $\min(W_J \d(w)) \le y \i \ast \d(x)$. \qed

\subsection{Admissible partitions of $\oG$}\label{subsec-admissible}
In order to generalize Theorem \ref{thm-closure-inter} to $\oG$, we will introduce the notion ``admissible partitions'' 
and discuss some of their properties. 

\begin{dfn}\label{dfn-admissible}
{\em A  partition of $\oG$ is said to be {\it admissible} if it is of the form
\begin{equation}\label{eq-one-parti}
\oG = \bigsqcup_{J \subset \Gamma} \bigsqcup_{\alpha \in \ca_J} X_{J, \alpha},
\end{equation}
where for each $J \subset \Gamma$ and $\alpha \in \ca_J$, $X_{J, \alpha} \subset Z_J$ and 
\[
{\rm Codim}_{Z_K} X_{K, \alpha'} \ge {\rm Codim}_{Z_J} X_{J, \alpha}
\]
for every $K \subset J$ and $X_{K, \alpha'} \subset 
\overline{X_\alpha} \cap Z_K$. 
An admissible partition is said to be {\it strongly admissible} if 
$\overline{X_{J, \alpha}} \cap Z_K \neq
\emptyset$ for every $K \subset J$ and $\alpha \in \ca_J$.
}
\end{dfn}

Note that the partition $\oG = \bigsqcup_{J \subset \Gamma} Z_J$ is strongly admissible.

\begin{prop}\label{pr-Z-O-admissible} The partitions 
\begin{align}\label{eq-4-part}
\oG &=\bigsqcup_{J \subset \G, (x, y) \in W^J \times W} [J, x, y]=\bigsqcup_{J \subset \G, (x, y) \in W^J \times W} 
[J, x, y]^{-, -} \\ 
\nonumber
&=\bigsqcup_{J \subset \G, (x, y) \in W^J \times W} [J, x, y]^{-, +}=\bigsqcup_{J \subset \G, w \in W^J} Z_{J, \d, w}
\end{align}
are strongly admissible.
\end{prop}

{\bf Proof.}
Let $K \subset J \subset \Gamma$ and $(x, y)\in W^J\times W$. If $(x', y') \in W^K \times W$ is  such that 
$[K, x', y'] \subset \overline{[J, x, y]}$, one knows
from  $\S$\ref {cl'} that there exists $u \in W_J$ such that $x' \ge x u$ and $y' \le y u$. Hence \begin{align*}
 {\rm Codim}_{Z_K} [K, x', y'] &=l(w_0)+l(x')-l(y') \ge l(w_0)+l(x u)-l(y u) \\
&=l(w_0) +l(x)+ l(u)-l(yu)\ge l(w_0) +l(x) -l(y)\\
&= {\rm Codim}_{Z_J} [J, x, y].
\end{align*} 
Regard$x$ as in $W^K$. By
$\S$\ref{cl'},  
$[K, x, y] \subset \overline{[J, x, y]} \cap Z_K$.
Thus the first partition in (\ref{eq-4-part}) is strongly admissible. 
The second and third partitions of $\oG$ in (\ref{eq-4-part}), being 
the
translations by $(w_0, 1)$ and by $(w_0, w_0)$ of the first one,  are thus also strongly admissible.

Consider now the partition of $\oG$ into the $G_\delta$-stable pieces. Let $K \subset J \subset \Gamma$ and $w \in W^J$. 
If  $w' \in W^K$ is such that
$Z_{K, \d, w'} \subset \overline{Z_{J, \d, w}}$, one knows from $\S$\ref{cl'}  that 
there exists $w_1 \in {\rm Min}(C_J(w))$ such that $w' \ge w_1$. Hence 
\[
{\rm Codim}_{Z_K}(Z_{K, \d, w'})=l(w') \ge l(w_1)=l(w)={\rm Codim}_{Z_J} (Z_{J, \d, w}).
\]
Regard $w$ as in $W^K$. By $\S$\ref{cl'}, $Z_{K, \d, w} \subset \overline{Z_{J, \d, w}} \cap Z_K$.
Thus the partition $\oG=\bigsqcup_{J \subset \G, w \in W^J} Z_{J, \d, w}$ is strongly admissible. 
 \qed 

\begin{nota-dfn}\label{dfn-compatible}
{\em Two admissible partitions
\begin{equation}\label{eq-two-parti}
\oG = \bigsqcup_{J \subset \Gamma} \bigsqcup_{\alpha \in \ca_J} X_{J, \alpha}
=\bigsqcup_{J \subset \Gamma} \bigsqcup_{\beta \in \cb_J} Y_{J, \beta}
\end{equation}
of $\oG$ are said to be {\it compatible} if for any $J \subset \Gamma$, $\alpha \in \ca_J$, and $\beta \in \cb_J$
with $X_{J, \alpha} \cap Y_{J, \beta} \neq \emptyset$, $X_{J, \alpha}$ and $Y_{J, \beta}$ intersect 
transversally in $Z_J$ and $X_{J, \alpha} \cap Y_{J, \beta}$ is irreducible.
For two such partitions of $\oG$,  and for $K \subset J \subset \Gamma$, $\alpha \in \ca_J$,
and $\beta \in \cb_J$, let
\begin{align*}
\ca_K^\alpha &= \{\alpha' \in \ca_K: \; X_{K, \,\alpha'} \subset \overline{X_{J,\, \alpha}},\;\;
{\rm Codim}_{Z_{K}} X_{K, \alpha'} = {\rm Codim}_{Z_J} X_{J, \alpha}\},\\
\cb_K^\beta &= \{\beta' \in \cb_K: \; Y_{K, \,\beta'} \subset \overline{Y_{J,\, \beta}},\;\;
{\rm Codim}_{Z_{K}} Y_{K, \beta'} = {\rm Codim}_{Z_J} Y_{J, \beta}\}.
\end{align*}}
\end{nota-dfn}

\begin{prop}\label{pr-comp}

1)  Any admissible partition of $\oG$ is compatible with the 
partition of $\oG$ into $G \times G$-orbits;

2)  The partitions
of $\oG$ into $\Gdel$-stable pieces and into $B^- \times B$-orbits are compatible;

3)  The partitions of $\oG$ into $B \times B$-orbits and into $B^- \times B^-$-orbits are compatible.
\end{prop}

{\bf Proof.} Directly from the definition and from Lemma \ref{lem-Richardson}. \qed

\

Recall \cite[Page 427]{Ht} that two irreducible subvarieties $X$ and $Y$ of a smooth 
irreducible variety $Z$ with $X \cap Y \neq \emptyset$ are said to
{\it intersect properly} in $Z$ if every irreducible component of 
$X \cap Y$ has dimension equal to $\dim X + \dim Y - \dim Z$.


\begin{thm}\label{thm-inter-admissible}
Let two compatible partitions of $\oG$ be given as in (\ref{eq-two-parti}). Then for any 
$J, K \subset \Gamma$ and $\alpha \in \ca_J$ and $\beta \in \cb_K$, if $\overline{X_{J, \alpha}}
\cap \overline{Y_{K, \beta}} \neq \emptyset$, then $\overline{X_{J, \alpha}}$ and 
$\overline{Y_{K, \beta}}$ intersect properly in $\overline{Z_{J \cup K}}$, and 
\[
\overline{X_{J, \alpha}}
\cap \overline{Y_{K, \beta}} = \bigcup_{(\alpha', \beta') \in \ci_{J \cap K}^{\alpha, \, \beta}}
\overline{X_{J\cap K, \, \alpha'} \cap Y_{J \cap K, \, \beta'}}
\]
is the decomposition of $\overline{X_{J, \alpha}}
\cap \overline{Y_{K, \beta}}$ into (distinct) irreducible components, where 
\[
\ci_{J \cap K}^{\alpha, \, \beta}=\{(\alpha', \beta') \in \ca_{J \cap K}^\alpha \times \cb_{J \cap K}^{\alpha, \beta}: \; 
X_{J\cap K, \, \alpha'} \cap Y_{J \cap K, \, \beta'} \neq \emptyset\}.
\]
In particular, $\ci_{J \cap K}^{\alpha, \, \beta}\neq \emptyset$.
 \end{thm}

{\bf Proof.} Let $J, K \subset \Gamma$, $\alpha \in \ca_J$, and $\beta \in \cb_K$ be such that $\overline{X_{J, \alpha}}
\cap \overline{Y_{K, \beta}} \neq \emptyset$. Regard both $\overline{X_{J, \alpha}}$ and 
$\overline{Y_{K, \beta}}$  as subvarieties of $\overline{Z_{J \cup K}}$. Since
$\overline{Z_{J \cup K}}$ is smooth and irreducible with
\[
\dim Z_{J \cup K} = \dim Z_J + \dim Z_K - \dim Z_{J \cap K}, 
\]
every
irreducible component of $\overline{X_{J, \alpha}}\cap \overline{Y_{K, \beta}}$ has dimension at least
\begin{align*}
l &= \dim X_{J, \alpha} + \dim Y_{K, \beta} - \dim Z_{J \cup K}\\
&=\dim Z_{J \cap K} - {\rm Codim}_{Z_J} X_{J, \,\alpha}
-{\rm Codim}_{Z_K} Y_{K, \,\beta}.
\end{align*}
On the other hand, 
\begin{equation}\label{eq-XY-inter}
\overline{X_{J, \alpha}}\cap \overline{Y_{K, \beta}} = \bigsqcup_{\begin{array}{c}I \subset J \cap K, \alpha' \in \ca_I, \beta' \in \cb_I \\
X_{I, \alpha'} \subset \overline{X_{J, \alpha}}, \, \, 
Y_{I, \beta'} \subset \overline{Y_{K, \beta}}\end{array}} X_{I, \alpha'} \cap Y_{I, \beta'}.
\end{equation}
For each non-empty intersection on the right hand side of (\ref{eq-XY-inter}), 
\begin{align}\label{eq-XY-dim}
\dim X_{I, \alpha'} \cap Y_{I, \beta'} & = \dim Z_I - {\rm Codim}_{Z_I} X_{I, \alpha'} -{\rm Codim}_{Z_I}  Y_{I, \beta'}\\
\nonumber& \le \dim Z_I - {\rm Codim}_{Z_J} X_{J, \alpha}-{\rm Codim}_{Z_K}  Y_{K, \beta}\\
\nonumber&\le \dim Z_{J \cap K} - {\rm Codim}_{Z_J} X_{J, \alpha}-{\rm Codim}_{Z_K}  Y_{K, \beta}\\
\nonumber &=l.
\end{align}
By Lemma \ref{le-general}, every irreducible component of $\overline{X_{J, \alpha}}\cap \overline{Y_{K, \beta}}$
has dimension $l$,  and the irreducible components are exactly as described in Theorem \ref{thm-inter-admissible}.
\qed

\
By taking the second admissible partition in Theorem \ref{thm-inter-admissible} to be the one by $G \times G$-orbits, we have
the following Corollary \ref{co-irr-1}.

\begin{cor}\label{co-irr-1}
Let a strongly admissible partition of $\oG$ be given as in (\ref{eq-one-parti}), and let $J \subset \Gamma$ and 
$\alpha \in \ca_J$. Then for any $K \subset \Gamma$, $\overline{X_{J, \alpha}} \cap \overline{Z_K} \neq \emptyset$ and  
$\overline{X_{J, \alpha}}$ and  $\overline{Z_K}$ intersect 
properly in $\overline{Z_{J \cup K}}$. Moreover,
\[
\overline{X_{J, \alpha}} \cap \overline{Z_K} = \bigcup_{\alpha' \in \ca_{J \cap K}^\alpha} \overline{X_{J \cap K, \, \alpha'}}
\]
is the decomposition of $\overline{X_{J, \alpha}} \cap \overline{Z_K}$ into (distinct) irreducible components.
\end{cor}

\begin{cor}\label{cor-irr-0}
Let two compatible partitions of $\oG$ be given as in (\ref{eq-two-parti}). Then for any $J \subset \Gamma$ and
$\alpha, \beta \in \ca_J$, 
$X_{J, \alpha} \cap Y_{J, \beta} \neq \emptyset$ if and only if $\overline{X_{J, \alpha}} \cap 
\overline{Y_{J, \beta}} \neq \emptyset$, and in this case,
\begin{equation}\label{eq-CC}
\overline{X_{J, \alpha}} \cap \overline{Y_{J, \beta}} = \overline{X_{J, \alpha} \cap Y_{J, \beta}}.
\end{equation}
In particular, 
\[
\oG = \bigsqcup_{J \subset \Gamma, \; \; 
X_{J, \alpha} \cap Y_{J, \beta} \neq \emptyset} X_{J, \alpha} \cap Y_{J, \beta}
\]
is an admissible partition of $\oG$.
\end{cor}

{\bf Proof.} Take $K = J$ in Theorem \ref{thm-inter-admissible}. If
$\overline{X_{J, \alpha}} \cap 
\overline{Y_{J, \beta}} \neq \emptyset$, then  $\cc_{K, \alpha, \beta}$ consists
of one element, namely, $(\alpha, \beta)$. Thus $X_{J, \alpha} \cap Y_{J, \beta} \neq \emptyset$ and
(\ref{eq-CC}) holds. The condition on codimensions in Definition \ref{dfn-admissible} follows
from (\ref{eq-XY-dim}) in the proof of Theorem \ref{thm-inter-admissible}.
\qed

\

\subsection{}
We now prove our second main result in this paper. Let
\begin{align*}
\cj &= \{(J, x, y, u, v): \;\;J \subset \Gamma, \,  (x, y), (u, v) \in W^J \times W,\\
& \hspace{1.3in} [J, x, y] \cap [J, u, v]^{-, -}\neq \emptyset\}\\
&=\{(J, x, y, u, v): \;\;J \subset \Gamma, \,  x, u \in W^J, y, v \in W,\\
& \hspace{1.3in} x \le u, \, v \le \max(yW_J)\},\\
\ck &= \{(J, w, x, y): \;\; J \subset \Gamma, \; (w, x, y) \in W^J \times W^J \times W, \\
& \hspace{1.3in}Z_{J, \delta, w} \cap [J, x, y]^{-, +}\neq \emptyset\}\\
&=\{(J, w, x, y): \;\; J \subset \Gamma, \; (w, x, y) \in W^J \times W^J \times W, \\
& \hspace{1.3in} \min (W_{ J} \, \delta(w)) \le  y^{-1} \ast \delta(x)\}.
\end{align*}

\begin{thm}\label{thm-BB0}
Let $J \subset \Gamma$ and $(x, y), (u, v) \in W^J \times W$. Then 
\[
[J, x, y] \cap [J, u, v]^{-, -} \neq \emptyset\;\;\; \mbox{iff} \;\;\;
\overline{[J, x, y]} \cap \overline{[J, u, v]^{-, -}} \neq \emptyset,
\]
 and in this case, 
\begin{equation}\label{xx} 
\overline{[J, x, y] \cap [J, u, v]^{-, -}} = \overline{[J, x, y]} \cap \overline{[J, u, v]^{-, -}}.
\end{equation}
In particular,
\begin{equation}\label{eq-JJ}
\oG = \bigsqcup_{(J, x, y, u, v) \in \cj} 
[J, x, y] \cap [J, u, v]^{-, -}
\end{equation}
is a strongly admissible partition of $\oG$.
\end{thm}

{\bf Proof.} Assume that  $\overline{[J, x, y]} \cap \overline{[J, u, v]^{-, -}} \neq \emptyset$. It follows from 
Corollary \ref{cor-irr-0} that $[J, x, y] \cap [J, u, v]^{-, -} \neq \emptyset$ and (\ref{xx}) holds. 

By Corollary \ref{cor-irr-0}, the partition (\ref{eq-JJ}) of $\oG$ is admissible. To show that
it is also strongly admissible, let $(J, x, y, u, v) \in \cj$ and let $K \subset J$.
By definition, there exists $z \in W_J$ such that $y z=\max(y W_J)$. Set $z'=\min(z W_K) \in W^K$. 
Then $x z', u z' \in W^K$ and 
\[
v z' \le \max(v W_J) \le \max(y W_J)= yz \le \max(y z' W_K).
\] 
Then $[K, x z', y z'] \cap [K, u z', v z']^{-, -} \neq \emptyset$. 
By $\S$\ref{cl'}, $[K, x z', y z'] \subset \overline{[J, x, y]}$ and 
$[K, u z', v z']^{-, -} \subset \overline{[J, u, v]^{-, -}}$. 
Therefore 
\[
\overline{[J, x, y] \cap [J, u, v]^{-, -}} \cap Z_K \supset [K, x z', y z'] \cap [K, u z', v z']^{-, -} \neq \emptyset.
\]
 This shows that the partition (\ref{eq-JJ}) is strongly admissible.
\qed

\begin{thm}\label{thm-inter-oG-2}
Let $J \subset \Gamma$, $w \in W^J$, and $(x, y) \in W^J \times W$. Then 
\[
Z_{J, \delta, w} \cap [J, x, y]^{-, +} \neq \emptyset
\;\;\; \mbox{iff} \;\;\; \overline{Z_{J, \delta, w}} \cap \overline{[J, x, y]^{-, +}}
\neq \emptyset,
\]
and in this case,
\begin{equation}\label{eq-OO}
\overline{Z_{J, \delta, w}} \cap \overline{[J, x, y]^{-, +}}=\overline{Z_{J, \delta, w} \cap [J, x, y]^{-, +}}.
\end{equation}
In particular,  
\begin{equation}\label{eq-oG-cap}
\oG = \bigsqcup_{(J, w, x, y) \in \ck} Z_{J, \delta, w} \cap [J, x, y]^{-, +}
\end{equation}
is a strongly admissible partition of $\oG$. 
\end{thm}

{\bf Proof.} Assume that  $\overline{Z_{J, \delta, w}} \cap \overline{[J, x, y]^{-, +}} \neq \emptyset$. It follows from 
Corollary \ref{cor-irr-0} that $Z_{J, \delta, w} \cap [J, x, y]^{-, +} \neq \emptyset$ and (\ref{eq-OO}) holds.

By Corollary \ref{cor-irr-0}, the partition (\ref{eq-oG-cap}) of $\oG$ is admissible. To show that
it is also strongly admissible, let $(J, w, x, y)\in \ck$ and let $K \subset J$.
By definition, there exists $z \in W_J$ such that $y z=\max(y W_J)$. Set $z'=\min(z W_K) \in W^K$. Then $x z' \in W^K$. 
Let $z = z' z^{\prime\prime}$ with $z^{\prime\prime} \in W_K$. Then 
$w_0^K \ast (y z') \i= w_0^K \ast (z^{\prime\prime} (yz)^{-1})= w_0^K \ast (yz)^{-1} = w_0^J \ast y \i.$ Thus
\begin{align*}
w_0^K \ast (y z') \i \ast \d(x z')&=w_0^J \ast y \i \ast \d(x z') 
\ge w_0^J \ast (y \i \ast \d(x))\\ &=\max(W_J (y \i \ast \d(x)).
\end{align*}
Since $\min(W_J \d(w)) \le y \i \ast \d(x)$, one has
\[
\d(w) \le \max(W_J (y \i \ast \d(x)) \le w_0^K \ast ((y z') \i \ast \d(x z')).
\]
Hence
$\min(W_K \d(w)) \le (y z') \i \ast \d(x z')$, and 
$$\overline{Z_{J, \d, w} \cap [J, x, y]^{-, +}} \cap Z_K \supset Z_{K, \d, w} \cap [K, x z', y z']^{-, +} \neq \emptyset. $$
This shows that the partition (\ref{eq-oG-cap}) of $\oG$ is strongly admissible.
\qed

\begin{rmk}\label{rem-another-check}
{\em Assume that $\overline{[J, x, y]} \cap \overline{[J, u, v]^{-, -}} \neq \emptyset$, we can also use 
Proposition \ref{prop-Z-O-inter-oG} to prove directly that $[J, x, y] \cap [J, u, v]^{-, -} \neq \emptyset$. 
Similarly, assume that $\overline{Z_{J, \delta, w}} \cap \overline{[J, x, y]^{-, +}}
\neq \emptyset$. One can use Proposition \ref{prop-Z-O-inter-oG} to prove directly that 
$Z_{J, \delta, w} \cap [J, x, y]^{-, +} \neq \emptyset$. 
We omit the details. 
}
\end{rmk}

Consider now the following four strongly admissible partitions of $\oG$:
\begin{align}\label{eq-4-parti}
\oG &=\bigsqcup_{J \subset \G, (x, y) \in W^J \times W} [J, x, y]^{-, +}=\bigsqcup_{J \subset \G, w \in W^J} Z_{J, \d, w}\\
\nonumber
&=\bigsqcup_{(J, x, y, u, v) \in \cj} [J, x, y] \cap [J, u, v]^{-, -}  =
\bigsqcup_{(J, x, y, u, v) \in \ck} Z_{J, \d, w} \cap [J, x, y]^{-, +}.
\end{align}

As a direct consequence of Corollary \ref{co-irr-1}, we have 
\begin{cor}\label{co-properly}
Let $J \subset \Gamma$ and let $X$ be any of the subvarieties of $Z_J$ appearing in either one of the fours partitions 
in (\ref{eq-4-parti}). Then for any $K \subset \Gamma$,
$\overline{X} \cap \overline{Z_K} \neq \emptyset$, and  $\overline{X}$ and $\overline{Z_K}$ intersect
properly in $\overline{Z_{J \cup K}}$.
\end{cor}

Corollary \ref{co-irr-1} also allows us to describe the irreducible components of the non-empty intersections in 
Corollary \ref{co-properly}.

\begin{cor}\label{co-ZO-ZK}
1) For any $J \subset \G$, $(x, y) \in W^J \times W$, and  $K \subset \G$, 
the irreducible components of $\overline{[J, x, y]^{-, +}} \cap \overline{Z_K}$ 
are precisely of the form $\overline{[J\cap K, x u, y u]^{-,+}}$, where $u \in W_J \cap W^{J \cap K}$ and $l(y u)=l(y)+l(u)$. 

2) For any $J \subset \Gamma$, $w \in W^J$, and  $K \subset \G$, the irreducible
components of $\overline{Z_{J, \delta, w}} \cap \overline{Z_K}$ are precisely of the form $\overline{Z_{J \cap K, \delta, w'}}$
with $w' \in W^{J \cap K} \cap {\rm Min}(C_J(w))$.

3) For any $(J, x, y, u, v) \in \cj$ and  $K \subset \Gamma$, the irreducible
components of the intersection $\overline{[J, x, y]
\cap [J, u, v]^{-, -}} \cap \overline{Z_K}$ are the non-empty
intersections of irreducible components of $\overline{[J, x, y]} \cap \overline{Z_K}$ and
the irreducible components of $\overline{[J, u, v]^{-, -} }\cap \overline{Z_K}$.

4) For any $(J, w, x, y) \in \ck$ and  $K \subset \Gamma$, the irreducible
components of $\overline{Z_{J, \delta, w} \cap [J, x, y]^{-, +}} \cap \overline{Z_K}$ are the non-empty
intersections of irreducible components of $\overline{Z_{J, \delta, w}} \cap \overline{Z_K}$ and
the irreducible components of $\overline{[J, x, y]^{-, +}} \cap \overline{Z_K}$.
\end{cor}

\begin{rmk}
{\em Corollary \ref{co-properly} and Corollary \ref{co-ZO-ZK}
in the case of the intersections $\overline{[J, x, y]^{-, +}} \cap \overline{Z_K}$ have also been obtained 
by M. Brion in
\cite{B} (using $B \times B^-$-orbits instead of $B^- \times B$-orbits). }
\end{rmk}

Applying Theorem \ref{thm-inter-admissible}, we have

\begin{cor}\label{cor-inter-ZJK}
1) If $J, K \subset \Gamma$ and $(x, y) \in W^J \times W, (u, v) \in W^K \times W$ are such that
$\overline{[J, x, y]} \cap \overline{[K, u, v]^{-, -}} \neq \emptyset$, then $\overline{[J, x, y]}$ and $ \overline{[K, u, v]^{-, -}}$ 
intersect properly in $\overline{Z_{J \cup K}}$, and the irreducible components of 
$\overline{[J, x, y]} \cap \overline{[K, u, v]^{-, -}} $ are the non-empty 
intersections of irreducible components of $\overline{[J, x, y]} \cap \overline{Z_{J \cap K}}$ and 
the irreducible components of $\overline{[K, u, v]^{-, -}} \cap \overline{Z_{J \cap K}}$.

2) If $J, K \subset \Gamma$ and $(w, x, y) \in W^J \times W^K \times W$ are such that  
$\overline{Z_{J, \delta, w}} \cap \overline{[K, x, y]^{-, +}} \neq \emptyset$, then  
$\overline{\Zw}$ and $\overline{[K, x, y]^{-, +}}$ intersect properly in $\overline{Z_{J \cup K}}$, and 
the irreducible components of 
$\overline{Z_{J, \delta, w}} \cap \overline{[K, x, y]^{-, +}}$ are the
non-empty intersections of irreducible components of $\overline{Z_{J, \delta, w}} \cap \overline{Z_{J \cap K}}$
and the irreducible components of $\overline{[K, x, y]^{-, +}} \cap \overline{Z_{J \cap K}}$.

\end{cor}

\section{A generalization of Deligne-Lusztig varieties}\label{sec-DL}

\subsection{} Let the ground field be an algebraically closed field in positive characteristic. 
Let $F: G \to G$ be a Frobenius map. We may choose a Borel subgroup $B$ and a maximal torus $H$ in 
such a way that $F(B)=B$ and $F(H)=H$. Then $F$ induces an automorphism on $W$ which we still 
denote by $F$. Set $G_F=\{(g, F(g)): g \in G\} \subset G \times G$. Let $\cc=(J, J', c, L)$ 
be an admissible quadruple as in $\S$\ref{sec-2.1}.
For $w \in W^J$, define $Z_{\cc, F, w}=G_F (B w, B)_\cdot R_\cc$. By \cite{H4}, 

1) $Z_{\cc, F, w}$ is a single $G_F$-orbit.

2) $Z_\cc=\bigsqcup_{w \in W^J} Z_{\cc, F, w}$.

3) $\overline{Z_{\cc, F, w}}=\bigsqcup Z_{\cc, F, w'}$, where $w'$ runs over elements in $W^J$ 
such that $\d \i(F(u)) w' u \i \le w$ for some $u \in W_J$. 

\subsection{} We now consider the intersection of a $G_\d$-stable piece and a $G_F$-orbit in $Z_\cc$. 
In the special case where $Z_\cc=G/B \times G/B$ and $\d$ is identity map, $Z_{\cc, \d, 1} \cap Z_{\cc, F, w'}$ 
are just the Deligne-Lusztig varieties \cite{DL}. 
It is also worth mentioning that in general Deligne-Lusztig varieties are not irreducible and not Frobenius split. See \cite{Ha}.  

\begin{prop}\label{pr-DL}
For any $w, w' \in W^J$, one has $Z_{\cc, \d, w} \cap Z_{\cc, F, w'} \neq \emptyset$, 
and $\overline{Z_{\cc, \d, w}} \cap \overline{Z_{\cc, F, w'}}=\overline{Z_{\cc, \d, w} \cap Z_{\cc, F, w'}}$. 
\end{prop}

{\bf Proof.} Notice that  $R_\cc \in Z_{\cc, \delta, 1} \cap Z_{\cc, F, 1}$. 
Hence $\overline{Z_{\cc, \d, w}} \cap \overline{Z_{\cc, F, w'}} \neq \emptyset$. 
Notice that ${\rm Lie}(G_F)=({\rm Lie}(G), 0)$. Thus ${\rm Lie}(G_\d)+{\rm Lie}(G_F)={\rm Lie}(G) \oplus {\rm Lie}(G)$. 
So $\overline{Z_{\cc, \d, w}}$ intersects transversally with $\overline{Z_{\cc, F, w'}}$. 
Therefore each irreducible component of $\overline{Z_{\cc, \d, w}} \cap \overline{Z_{\cc, F, w'}}$ 
is of dimension $\dim\overline{Z_{\cc, \d, w}}+\dim\overline{Z_{\cc, F, w'}}-\dim Z_\cc$. 
On the other hand, $\overline{Z_{\cc, \d, w}} \cap \overline{Z_{\cc, F, w'}}$ is the union of 
subvarieties $Z_{\cc, \d, x} \cap Z_{\cc, F, y}$, where $x$ and $y$ run over elements in $W^J$ 
such that $Z_{\cc, \d, x} \subset \overline{Z_{\cc, \d, w}}$ and 
$Z_{\cc, F, y} \subset \overline{Z_{\cc, F, w'}}$. For such $x$ and $y$, 
\begin{align*} 
\dim(Z_{\cc, \d, x} \cap Z_{\cc, F, y})&=\dim Z_{\cc, \d, x}+\dim Z_{\cc, F, y}-\dim Z_\cc \\ 
& \le \dim\overline{Z_{\cc, \d, w}}+\dim\overline{Z_{\cc, F, w'}}-\dim Z_\cc 
\end{align*} 
with equality holds only when $\dim Z_{\cc, \d, x}=\dim\overline{Z_{\cc, \d, w}}$ and 
$\dim Z_{\cc, F, y}=\dim\overline{Z_{\cc, F, w'}}$, i.e., $x=w$ and $y=w'$. 
Therefore the irreducible components of $\overline{Z_{\cc, \d, w}} \cap \overline{Z_{\cc, F, w'}}$ 
are precisely the closure $\overline{Y}$, where $Y$ is an irreducible component of $Z_{\cc, \d, w} \cap Z_{\cc, F, w'}$. 
In particular, $Z_{\cc, \d, w} \cap Z_{\cc, F, w'} \neq \emptyset$ and 
$\overline{Z_{\cc, \d, w}} \cap \overline{Z_{\cc, F, w'}}=\overline{Z_{\cc, \d, w} \cap Z_{\cc, F, w'}}$. \qed

\begin{rmk}
{\em The Proposition \ref{pr-DL} can also be generalized to $\oG$. We omit the details.}
\end{rmk}




\section{Appendix}

Recall that $W$ is the Weyl group of $G$. We now prove some properties of the operations $\ast, \tril, \trir$ 
on $W$ as defined in Section \ref{ast}. In fact, many properties also holds for arbitrary Coxeter groups. See \cite{H5}. 

\begin{lem}\label{lem-he}\cite[Lemma 3.3]{H2}
For any $x, y \in W$,

1) $x \ast y \in W$ is the unique maximal element in the set  $\{u y: u \le x\}$
as well as in the set $\{x v; v \le y\}$. Moreover, $x \ast y = x_1 y = x y_1$, where
$x_1 \le x, y_1 \le y$ and $l(x \ast y) = l(x_1) + l(y) = l(x) + l(y_1)$;

2) $x \trir y \in W$ is the unique minimal element in the set $\{u y: u \le x\}$, and
$x \trir y = x_1 y$ with $x_1 \le x$ and $l(x \trir y) = l(y) - l(x_1)$;

3)  $x \tril y \in W$ is the unique minimal element in the set $\{x v; v \le y\}$, and
$x \tril y = x y_1$ with $y_1 \le y$ and $l(x \tril y) = l(x) - l(y_1)$.
\end{lem}

\begin{lem}\label{2}
Let $x, x', y, y' \in W$. If $x \le x'$ and $y \le y'$, then 
\[
x \ast y \le x' \ast y', \hs x' \trir y \le x \trir y', \hs \mbox{and} \hs x \tril y' \le x' \tril y.
\]
\end{lem}

{\bf Proof.} It follows from 
\[
B (x \ast y) B \subset BxByB \subset \overline{Bx'B} \; \overline{B y' B}
\subset \overline{B x' B y' B}=\overline{B (x' \ast y') B},
\]
that   $x \ast y \le x' \ast y'$. Similarly, since
\[
B (x \trir y') B^- \subset BxBy'B^- \subset \overline{Bx'B} \; 
\overline{ByB^-} \subset \overline{ B x' B y B^-} =
\overline{B (x' \trir y) B^-},
\] 
one has $x' \trir y \le x \trir y'$. 
Similarly, $x \tril y' \le x' \tril y$.
\qed

\begin{lem}\label{lem-facts1} For any $x, y, z \in W$, 

1) $x \trir y =(x \ast (yw_0))w_0$ and $x \tril y =  w_0 ( (w_0x) \ast y)$;

2) $(x \tril y)^{-1} = y^{-1} \trir x^{-1}$ and $(x \ast y)^{-1}=y \i \ast x \i$;

3) $x \trir y \le z$ if and only if $y \le x^{-1} \ast z$;

4) $y \tril x \le z$ if and only if $y \le z \ast x^{-1}$;

5) $(x \trir y) \tril z = x \trir (y \tril z)$.
\end{lem}

{\bf Proof.} 1) Since 
\begin{align*} \overline{ B (x \trir y) B^-} 
&=\overline{ B x B y B^-}=\overline{B x B y w_0 B w_0}=\overline{B (x \ast (y w_0)) B w_0} \\ 
&=\overline{B \bigl(x \ast (y w_0)\bigr) w_0 B^-},
\end{align*} 
one has $x \trir y=x \trir (y w_0) w_0$. 
Similarly, $x \tril y =  w_0 ( (w_0x) \ast y)$.

2) Let $\tau$ be the inverse map of $G$. Then $(x \tril y)^{-1} = y^{-1} \trir x^{-1}$
follows by applying $\tau$ to $\overline{B^- (x \tril y) B} = \overline{B^- x B y B}$.
Similarly, $(x \ast y)^{-1}=y \i \ast x \i$.

3) Since $y \in \{ u (x \trir y): u \le x \i\}$, $y \le x \i \ast (x \trir y)$.
If $x \trir y \le z$, then $y \le x \i \ast z$ by Lemma \ref{2}.
Similarly, $z \in \{u (x \i \ast z): u \le x\}$, so $x \trir (x \i \ast z) \le z$.
If $y \le x \i \ast z$, then by Lemma \ref{2}, $x \trir y \le x \trir (x^{-1} \ast z) \le z$.

Part 4) can be proved in the same way as part 3).

5) Since
\begin{align*}
B ((x \trir y) \tril z) B^- & \subset B (x \trir y) B^-z B^- \subset  
B x B y B^- z B^- \\
& \subset \overline{BxB} \; \overline{B y B^- z B^-}  = 
\overline{BxB} \; \overline{B (y \tril z) B^-}\\
& \subset \overline{BxB (y \tril z) B^-} =  \overline{B( x \trir (y \tril z)) B^-},
\end{align*}
one sees that $(x \trir y) \tril z \ge x \trir (y \tril z)$. Similarly, one shows that
$(x \trir y) \tril z \le x \trir (y \tril z)$. Thus $(x \trir y) \tril z = x \trir (y \tril z)$.
\qed




\begin{lem}\label{lem-tri-min}
For $J, J' \subset \Gamma$, $x \in W$,  $y \in W^J$ and $z \in {}^{J}W$, one has
$x \trir y \in W^J$, $z \tril x \in {}^{J}W$, and
\[
w_0^{J'} \trir x \tril w_0^J =\min(W_{J'} \, x W_J), \hs w_0^{J'} \ast x \ast w_0^J =\max(W_{J'} \, x W_J).
\]
\end{lem}

{\bf Proof.} By Lemma \ref{lem-he}, $x \trir y = x_1 y$ with $x_1 \le x$ and
$l(x \trir y) = l(y) - l(x_1)$. By \cite[3.5]{H2}, $x_1 y \in W^J$. Similarly one
has $z \tril x \in {}^{J}W$. It is clear from the definitions that 
\[
w_0^{J'} \trir x = \min (W_{J'} x) \in {}^{J'}\!W \hs \mbox{and} \hs x \tril w_0^J = \min (x W_J) \in W^J.
\]
Thus $w_0^{J'} \trir x \tril w_0^J \in {}^{J'}\!W^J$. By Lemma \ref{lem-he}, $w_0^{J'} \trir x \tril w_0^J
\in W_{J'} \, x W_J$. Thus $w_0^{J'} \trir x \tril w_0^J =\min(W_{J'} \, x W_J).$
Similarly, $w_0^{J'} \ast x \ast w_0^J =\max(W_{J'} \, x W_J)$.
\qed

\

Combining Lemma \ref{lem-tri-min} with Lemma \ref{lem-facts1} 4), we have the following consequence.

\begin{lem}\label{lem-x} For any $J \subset \Gamma$ and $x, y \in W$, $x \le \max(y W_J)$ if and only if $\min(x W_J) \le y$.
\end{lem}

\

The following Lemma \ref{lem-deodhar} can be found in \cite{De, KL}. 

\begin{lem}\label{lem-deodhar} For $x, y \in W$, the following conditions are equivalent:

1) $B x B \subset B y B$;

2) $B^- y B \subset B^- x B$;

3) $(B^- x B) \cap (B y B) \neq \emptyset$;

4) $\overline{B^- x B} \cap \overline{B y B} \neq \emptyset$;

5) $x \le y$.
\end{lem}

\

The following result will be used several times in our paper. 

\begin{lem}\label{lem-facts3}
For $x, y, u, v \in W$, the following conditions are equivalent:

1) $(BxByB) \cap (B^-uBvB) \neq \emptyset$;

2) $\overline{BxByB} \cap (B^-uBvB) \neq \emptyset$;

3) $(BxByB) \cap \overline{B^-uBvB} \neq \emptyset$;

4) $\overline{BxByB} \cap \overline{B^-uBvB} \neq \emptyset$;

5) $u \tril v  \le x \ast y$.

6) $u \le x \ast y \ast v \i$. 
\end{lem}

{\bf Proof.} Clearly 1) implies 2) and 3), 2) or 3) implies 4),
4) implies 5) by Lemma \ref{lem-deodhar} and 5) is equivalent to 6) by Lemma \ref{lem-facts1}. It suffices to show that 5) implies 1).

Suppose that $u \tril v \le x \ast y$. Then $(B (x \ast y) B) \cap (B^- (u \tril v) B) \neq \emptyset$
by Lemma \ref{lem-deodhar}. 
Since $B (x \ast y) B \subset B x B y B$ and $B^- (u \tril v) B \subset B^- u B v B$, 
we have  $(B x B y B) \cap (B^- u B v B) \neq \emptyset$. Hence 5) implies 1).
\qed

\section*{Acknowledgments} We thank J. Starr and J. F. Thomsen for some helpful discussions. The first author is partially supported by HKRGC grants 601409 and DAG08/09.SC03. The second author is partially supported by HKRGC grants 703405 and 703707.


\begin{thebibliography}{99}
\bibitem{B} M. Brion, {\em The behavior at infinity of the Bruhat decomposition}, Comment. Math. Helv. {\bf 73}
(1998), 137 - 174.


\bibitem{DP} C. De Concini and C. Procesi, {\em{Complete symmetric
varieties,}} in Invariant theory (Montecatini, 1982), Lecture Notes
in Math. {\bf{996}} (1983), 1--44.

\bibitem{DS} C. De Concini and T. Springer, {\em{Compactification
of symmetric spaces,}} Transformation Groups {\bf{4}} (1999), 273-300.

\bibitem{De} V. Deodhar, {\em On some geometric aspects of Bruhat orderings, I. A finer decomposition of Bruhat cells}, Invent. Math. {\bf 79} (1985), 499 - 511.

\bibitem{DL}
P. Deligne and G. Lusztig, {\em Representations of reductive groups over finite fields}, Ann. of Math. (2) 103 (1976), no. 1, 103--161.

\bibitem{e-l:cplx}
S. Evens and  J.-H. Lu, {\em On the variety of Lagrangian subalgebras, II}, 
Ann. Scient. Ecol. Norm. Sup., {\bf 39} (2) (2006), 347 - 379.

\bibitem{e-l:grothendieck}
S. Evens and  J.-H. Lu, {\em Poisson geometry of the Grothendieck resolution of a 
complex semisimple group}, Moscow Mathematical Journal {\bf 7} (4) (special volume 
in honor of V. Ginzburg's 50'th birthday) (2007),  613 - 642. 

\bibitem{fz-double-positive}
S. Fomin and A. Zelevinsky,  {\em Double Bruhat cells and total positivity},
{\em J. Amer. Math. Soc.} {\bf 12} (1999), no. 2, 335--380.

\bibitem{Ha}
S. H. Hansen, {\em Canonical bundles of Deligne-Lusztig varieties}, Manuscripta Math. 98 (1999), 363--375.

\bibitem{Hr} J. Harris, {\em Algebraic geometry}, Graduate Texts in Mathematics, {\bf 133} Springer-Verlag, 1992.

\bibitem{Ht} R. Hartshorne, {\em Algebraic geometry}, Graduate Texts in Mathematics, {\bf 52}, Springer-Verlag, 1977.


\bibitem{H2}
X. He, {\em The $G$-stable pieces of the wonderful compactification}, Tran. of AMS {\bf 359} (7) (2007), 3005 - 3024.

\bibitem{H3}
X. He, {\em Minimal length elements in some double cosets of Coxeter groups}, Adv. Math. {\bf 215} (2007), 469 - 503.

\bibitem{H4}
X. He, {\em G-stable pieces and partial flag varieties}, Representation Theory, Contemp. Math., vol. 478, Amer. Math. Soc., Providence, RI, 2009, pp. 61Ð70. 

\bibitem{H5}
X. He, {\em A subalgebra of $0$-Hecke algebra}, to appear in J. Algebra, 2009. 


\bibitem{HT2} X. He and J. F. Thomsen,
{\em Geometry of $B \times B$-orbit closures in equivariant embeddings}, Adv. Math. {\bf 216} (2) (2007), 626--646.

\bibitem{HT3} X. He and J. F. Thomsen, {\em Frobenius splitting and geometry of $G$-Schubert varieties}, {\em Adv. Math.} {\bf 219} (5) (2008), 1469 - 1512.

\bibitem{KL}
D. Kazhdan and G. Lusztig, {\em Schubert varieties and Poincar\'e duality}, Proc. Symp. Pure Math. {\bf 36}, 
Amer. Math. Soc. (1980), pp. 185-203.

\bibitem{KZ}
M. Kogan  and A. Zelevinsky, {\em On symplectic leaves and integrable systems 
in standard complex semisimple Poisson-Lie groups}, Int. Math. Res. Not
(32) (2002), 1685--1702.

\bibitem{LY}
J. Lu and M. Yakimov, {\em Partitions of the wonderful group compactification}, Trans. Groups, 
{\bf 12} (4) (2007), 695 - 723.

\bibitem{LY2}
J. Lu and M. Yakimov, {\em Group orbits and regular partitions of Poisson
manifolds}, Comm. Math. Phys., {\bf 283}(3), 729 - 748 (2008). 

\bibitem{L1} G. Lusztig, {\em{Parabolic character sheaves I}},
Moscow Math J. {\bf{4}} (2004), 153-179.
       
\bibitem{L2} G. Lusztig, {\em{Parabolic character sheaves II}},
Moscow Math J. {\bf{4}} (2004), 869-896.


\bibitem{R}
R. W. Richardson, {\em Intersections of double cosets in algebraic
groups}, Indagationes Mathematicae, Volume 3, Issue 1, (1992),
69-77.

\bibitem{RS}
R. W. Richardson and T. A. Springer, {\em The Bruhat order on symmetric spaces},
Geom. Dedic.  {\bf 35} (1990), 389 - 436.

\bibitem{Ro}
R. Rouquier, {\em Weyl groups, affine {W}eyl groups and reflection groups}, Representations of reductive groups, 21--39, 
Publ. Newton Inst., Cambridge Univ. Press, Cambridge, 1998. 

\bibitem{Sp} T. A. Springer, {\em{Intersection cohomology of 
$B\times B$-orbit closures in group compactifications,}} with an appendix by 
W. van der Kallen. Special issue in celebration of Claudio Procesi's 
60th birthday, J. Algebra {\bf{258}} (2002), 71--111. 

\bibitem{Sp2}
T. A. Springer, {\em An extension of Bruhat's Lemma}, J. Alg. 313 (2007), 417--427.


\end{thebibliography}
\end{document}